\documentclass[11pt,a4paper,leqno]{article}
\usepackage{amsmath}
\usepackage{amssymb}
\usepackage{amsthm}
\usepackage{enumerate}
\usepackage{txfonts,url}
\usepackage[text={135mm,185mm},centering]{geometry}
\newcommand\tov[1][v]{\overset{#1}{\to}}
\newcommand\notov{\overset{v}{\nrightarrow}}
\newcommand\FN[1][q]{\ensuremath{F_v(a_1,\dots,a_r;#1)}}
\newcommand\cl{\operatorname{cl}}
\newcommand\ol{\overline}
\newcommand\MM{\mathcal{M}}

\numberwithin{equation}{section}
\numberwithin{table}{section}

\newtheorem{thm}{Theorem}[section]
\newtheorem{lemma}{Lemma}[section]
\newtheorem{cor}{Corollary}[section]
\theoremstyle{definition}
\newtheorem{definition}{Definition}[section]
\newtheorem{rem}{Remark}[section]

\begin{document}
\title{On the Vertex Folkman Numbers $F_v(2,\dots,2;q)$%
\thanks{Supported by the Scientific Research Fund of
St. Kl. Ohridski Sofia University under contract 90/2008.}}
\author{Nedyalko Nenov}
\maketitle

\begin{abstract}
For a graph $G$ the symbol $G\tov(a_1,\dots,a_r)$ means
that in every $r$-coloring of the vertices of $G$
for some $i\in\{1,\dots,r\}$ there exists a monochromatic
$a_i$-clique of color $i$.
The vertex Folkman numbers
\[
\FN=\min\{|V(G)|:G\tov(a_1,\dots,a_r)\text{ and }K_q\nsubseteqq G\}
\]
are considered. In this paper we shall compute the Folkman numbers
$F_v(\underbrace{2,\dots,2}_r;r-k+1)$ when $k\le 12$ and $r$ is sufficiently
large. We prove also new bounds for some vertex and edge Folkman numbers.

\textbf{2000 Mathematics Subject Classification:} 05C55

\textbf{Key words:} Folkman numbers, vertex coloring, edge coloring
\end{abstract}

\section{Introduction}\label{s:1}

We consider only finite, non-oriented graphs without loops and multiple edges.
The vertex set and the edge set of a graph $G$ will be denoted by $V(G)$ and
$E(G)$, respectively. A graph $G$ is said to be an \emph{empty graph} if
$V(G)=\emptyset$. We call \emph{a $p$-clique} of a graph $G$ a set of $p$
pairwise adjacent vertices. The largest integer $p$ such that the graph $G$
contains a $p$-clique is denoted by $\cl(G)$. A set of vertices of a graph
is said to be \emph{independent} if every two of them are not adjacent.
We shall also use the following notations:

$\ol G$ is the complement of $G$;

$\alpha(G)$ is the vertex independence number of $G$,
i.e., $\alpha(G)=\cl(\ol G)$;

$\chi(G)$ is the chromatic number of $G$;

$f(G)=\chi(G)-\cl(G)$;

$K_n$ is the complete graph on $n$ vertices;

$C_n$ is the simple cycle on $n$ vertices;

$M(x,y)=\{G:|V(G)|<\chi(G)+2f(G)-x\text{ and }f(G)\le y\}$.

The graph $G$ is a \emph{$(p,q)$-graph} if $\cl(G)<p$ and $\alpha(G)<q$.
The \emph{Ramsey number} $R(p,q)$ is the smallest natural $n$ such that
every graph $G$ with $|V(G)|\ge n$ is not a $(p,q)$-graph. An exposition
of the results on the Ramsey numbers is given in \cite{Rad06}. We shall
need Table~\ref{t:1.1} of the known Ramsey numbers $R(p,3)$ (see \cite{Rad06}).

\begin{table}[h]
\centering
\caption{Ramsey numbers $R(p,3)$}\label{t:1.1}
\begin{tabular}{*{10}{|c}|}
\hline
$p$&3&4&5&6&7&8&9&10&11\\
\hline
$R(p,3)$&6&9&14&18&23&28&36&40--43&46--51\\
\hline
\end{tabular}
\end{table}

Let $G_1$ and $G_2$ be two graphs without common vertices. We denote
by $G_1+G_2$ the graph $G$ for which $V(G)=V(G_1)\cup V(G_2)$ and
$E(G)=E(G_1)\cup E(G_2)\cup E'$ where $E'=\{[x,y], x\in V(G_1),
y\in V(G_2)\}$. A graph $G$ is \emph{separable} if $G=G_1+G_2$,
where $V(G_i)=\emptyset$, $i=1,2$.

\begin{definition}\label{d:1.1}
Let $\MM\ne\emptyset$ be a set of graphs. We say that $G_0\in\MM$ is
a \emph{minimal graph} in $\MM$ if $|V(G_0)|=\min\{|V(G)|:G\in\MM\}$.
\end{definition}

\begin{definition}\label{d:1.2}
Let $a_1,\dots,a_r$ be positive integers. The symbol $G\tov(a_1,\dots,a_r)$
means that in every $r$-coloring
\[
V(G)=V_1\cup\dots\cup V_r,
\qquad
V_i\cap V_j=\emptyset,
\quad
i\ne j
\]
of the vertices of $G$ for some $i\in\{1,\dots,r\}$ there exists
a monochromatic $a_i$-clique $Q$ of color $i$, that is $Q\subseteq V_i$.

Define
\begin{align*}
H_v(a_1,\dots,a_r;q)&=\{G\tov(a_1,\dots,a_r)\text{ and }\cl(G)<q\},\\
F_v(a_1,\dots,a_r;q)&=\min\{|V(G)|:G\in H_v(a_1,\dots,a_r;q)\}.
\end{align*}
\end{definition}

It is clear that $G\tov(a_1,\dots,a_r)$ implies $\cl(G)\ge\max\{a_1,\dots,a_r\}$.
Folkman proved in \cite{Folkman} that there exists a graph $G$ such that
$G\tov(a_1,\dots,a_r)$ and $\cl(G)=\max\{a_1,\dots,a_r\}$. Therefore,
\begin{equation}\label{1.1}
F_v(a_1,\dots,a_r;q)\text{ exists}\iff q>\max\{a_1,\dots,a_r\}.
\end{equation}

The numbers $F_v(a_1,\dots,a_r;q)$ are called \emph{vertex Folkman numbers}.
If $a_1,\dots,a_r$ are positive integers, $r\ge 2$ and $a_i=1$ then it is
easy to see that
\[
G\tov(a_1,\dots,a_{i-1},a_i,a_{i+1},\dots,a_r)
\Rightarrow
G\tov(a_1,\dots,a_{i-1},a_{i+1},\dots,a_r).
\]
Hence, for $a_i=1$
\[
F_v(a_1,\dots,a_r;q)=F_v(a_1,\dots,a_{i-1},a_{i+1},\dots,a_r;q).
\]
Thus, it is enough to consider just such numbers $F_v(a_1,\dots,a_r;q)$
for which $a_i\ge 2$. In this paper we consider the Folkman numbers
$F_v(2.\dots,2;q)$. Set
\[
(\underbrace{2,\dots,2}_r)=(2_r)
\text{ and }
F_v(\underbrace{2,\dots,2}_r;q)=F_v(2_r;q).
\]
By~\eqref{1.1}
\begin{equation}\label{1.2}
F_v(2_r;q)\text{ exists}\iff q\ge 3.
\end{equation}
It is clear that
\begin{equation}\label{1.3}
G\tov(2_r)\iff\chi(G)\ge r+1.
\end{equation}
Since $K_{r+1}\tov(2_r)$ and $K_r\overset{v}{\not\to}(2_r)$ we have
\begin{equation}\label{1.4}
F_v(2_r;q)=r+1
\text{ if }
q\ge r+2.
\end{equation}
According to~\eqref{1.4} it is enough to consider just such numbers
$F_v(2_r;r-k+1)$ for which $k\ge-1$. In this paper we shall prove
the following results.

\begin{thm}\label{th:1.1}
Let $r$ and $k$ be integers such that $-1\le k\le 5$ and $r\ge k+2$. Then
\begin{enumerate}[\indent\rm(a)]
\item
$F_v(2_r;r-k+1)\ge r+2k+3$;
\item
$F_v(2_r;r-k+1)=r+2k+3$ if $k\in\{0,2,3,4,5\}$ and $r\ge 2k+2$ or
$k\in\{-1,1\}$ and $r\ge 2k+3$.
\end{enumerate}
\end{thm}

\begin{thm}\label{th:1.2}
Let $r\ge 8$ be a natural number. Then
\begin{enumerate}[\indent\rm(a)]
\item
$F_v(2_r;r-5)\ge r+14$ and $F_v(2_r;r-5)=r+14$ if and only if $r\ge 13$;
\item
$F_v(2_r;r-6)\ge r+16$ if $r\ge 9$ and
$F_v(2_r;r-6)=r+16$ if $r\ge 15$;
\item
$F_v(2_r;r-7)\ge r+17$, $r\ge 10$ and
$F_v(2_r;r-7)=r+17$ if and only if $r\ge 16$;
\item
$F_v(2_r;r-8)\ge r+18$, $r\ge 11$ and
$F_v(2_r;r-8)=r+18$ if and only if $r\ge 17$;
\item
$F_v(2_r;r-9)\ge r+20$, $r\ge 12$ and
$F_v(2_r;r-9)=r+20$ if $r\ge 19$.
\end{enumerate}
\end{thm}

\begin{thm}\label{th:1.3}
Let $r\ge 13$ be a natural number. Then
\begin{enumerate}[\indent\rm(a)]
\item
$F_v(2_r;r-10)\ge r+21$ and $F_v(2_r;r-10)=r+21$ if $R(10,3)>41$
and $r\ge 20$;
\item
If $R(10,3)\le 41$ we have $F_v(2_r;r-10)\ge r+22$ and
$F_v(2_r;r-10)=r+22$ if $r\ge 21$.
\end{enumerate}
\end{thm}

\begin{thm}\label{th:1.4}
Let $r$ and $k$ be natural numbers such that $r\ge k+2$ and $k\ge 12$. Then
\begin{enumerate}[\indent\rm(a)]
\item
$F_v(2_r;r-k+1)\ge r+k+11$;
\item
If $k=12$ and $r\ge 22$ then $F_v(2_r;r-11)=r+23$.
\end{enumerate}
\end{thm}

\begin{rem}\label{rem:1.1}
By~\eqref{1.2} the number $F_v(2_r;r-k+1)$ exists if and only if $r\ge k+2$.
Thus, the inequality $r\ge k+2$ in the statements of these Theorems is
necessary.
\end{rem}

\begin{rem}\label{rem:1.2}
The case $k=0$ of Theorem~\ref{th:1.1} was proved by Dirac in \cite{Dirac}.
It was also proved in~\cite{Dirac} that the graph $K_{r-2}+C_5$, $r\ge 2$
is the only minimal graph in $H_v(2_r;r+1)$. The cases $k=1$ and $k=2$ of
Theorem~\ref{th:1.1} were proved in~\cite{Nen83}. It was also proved
in~\cite{Nen83} that $K_{r-5}+C_5+C_5$, $r\ge 5$ is the only minimal graph
in $H_v(2_r;r)$ (see also~\cite{Nen03}). The case $k=3$ was proved
in~\cite{Nen81a}. We gave new proofs of the cases $k=2$ and $k=3$ of
Theorem~\ref{th:1.1} in~\cite{NenAnn}.
\end{rem}

The method we use here does not allow us to compute the numbers
$F_r(2_r;r-k+1)$ when $r<2k+2$ and $1\le k\le 5$. We know only
the following numbers of this kind:
\begin{align*}
F_v(2_3;3)&=11,&&\text{\cite{Chvatal} and \cite{Myc};}\\
F_v(2_4;3)&=22,&&\text{\cite{JenRoy};}\\
F_v(2_r;4)&=11,&&\text{\cite{Nen84} (see also \cite{Nen98}).}
\end{align*}

We know about number $F_4(2_5;4)$ only that $12\le F_v(2_5;4)\le 16$
(see \cite{NenAnn}).

\begin{rem}\label{rem:1.3}
If $k\ge 2$ then there is more than one minimal graph in $H_v(2_r;r-1)$.
For example, if $r\ge 8$ the graph $K_{r-8}+C_5+C_5+C_5$ is also minimal
in $H_v(2_r;r-1)$ besides the minimal graph from the proof of
Theorem~\ref{th:1.1}.
\end{rem}

\begin{rem}\label{rem:1.4}
Luczak et~al. \cite{LucRucUrb} proved the inequality
\begin{equation}\label{1.5}
F_v(2_r;r-k+1)\le r+2k+3
\text{ if }
r\ge 3k+2.
\end{equation}
They also proved that~\eqref{1.5} is strict when $k$ is very large
(see~\cite{LucRucUrb}). It can be easily seen from Theorem~\ref{th:1.1}
and Theorem~\ref{th:1.2}~(a) that $k=6$ is the smallest value of $k$
for which the inequality~\eqref{1.5} is strict.
\end{rem}

\section{Auxiliary Results}\label{s:2}

The following lemmas are used to prove the main results.

\begin{lemma}\label{l:2.1}
Let $q\ge 4$ be an integer and $G$ be a minimal graph
(see Definition~\ref{d:1.1}) in $H_v(2_r;q-1)$. Then
\[
F_v(2_r;q-1)\ge F_v(2_r;q)+\alpha(G)-1.
\]
\end{lemma}

\begin{proof}
Let $A\subseteq V(G)$ be an independent set of vertices of $G$ such that
$|A|=\alpha(G)$. Consider the graph $G'=K_1+(G-A)$. By~\eqref{1.3},
$\chi(G)\ge r+1$. Since $A$ is an independent vertex set it follows that
$\chi(G-A)\ge r$ and $\chi(G')\ge r+1$. By~\eqref{1.3}, $G'\tov(2_r)$.
Since $\cl(G)\le q-2$ we have $\cl(G')\le q-1$. Hence, $G'\in H_v(2_r;q)$
and
\[
F_v(2_r;q)\le|V(G')|=|V(G)|-\alpha(G)+1.
\]
Lemma~\ref{l:2.1} follows from this inequality because $|V(G)|=F_v(2_r;q-1)$.
\end{proof}

\begin{cor}\label{cor:2.1}
Let $q$ and $r$ be integers such that $4\le q<r+3$. Then
\begin{enumerate}[\indent\rm(a)]
\item
$F_v(2_r;q-1)\ge F_v(2_r,q)+1$;
\item
If $F_v(2_r;q)+1\ge R(q-1,3)$ then the inequality \textup{(a)} is strict.
\end{enumerate}
\end{cor}

\begin{proof}
Let $G$ be a minimal graph in $H_v(2_r;q-1)$. By~\eqref{1.3}, $\chi(G)\ge r+1$.
Since $\cl(G)\le q-2$ and $q<r+3$ we have
\[
\cl(G)<r+1\le\chi(G).
\]
Thus, $\alpha(G)\ge 2$ and inequality~(a) follows from Lemma~\ref{l:2.1}.

Let $F_v(2_r;q)+1\ge R(q-1,3)$. Then we see from~(a) that
\[
|V(G)|=F_v(2_r;q-1)\ge R(q-1,3).
\]
Since $\cl(G)<q-1$, this inequality implies $\alpha(G)\ge 3$.
From Lemma~\ref{l:2.1} we obtain
\[
F_v(2_r;q-1)\ge F_v(2_r;q)+2.
\]
The Corollary~\ref{cor:2.1} is proved.
\end{proof}

A graph $G$ is said to be \emph{$k$-chromatic} if $\chi(G)=k$.
A graph $G$ is defined to be \emph{vertex-critical chromatic}
if $\chi(G-v)<\chi(G)$ for all $v\in V(G)$.

\begin{lemma}\label{l:2.2}
Let $q\ge 3$ be an integer and let $G$ be a minimal graph in $H_v(2_r;q)$.
Then
\begin{enumerate}[\indent\rm(a)]
\item
$G$ is a vertex-critical $(r+1)$-chromatic graph;
\item
If $q<r+3$ then $\cl(G)=q-1$.
\end{enumerate}
\end{lemma}

\begin{proof}
\emph{Proof of} ($a$).
By~\eqref{1.3}, $\chi(G)\ge r+1$. Assume that (a) is wrong.
Then there exists $v\in V(G)$ such that $\chi(G-v)\ge r+1$. Thus,
according to~\eqref{1.3}, $G-v\in H_v(2_r;q)$. This contradicts
the minimality of $G$ in $H_v(2_r;q)$.

\emph{Proof of} ($b$).
Assume that (b) is wrong. Then, since $\cl(G)\le q-1$ we have
$\cl(G)\le q-2$. Thus, $G\in H_v(2_r;q-1)$. Hence $H_v(2_r;q-1)
\ne\emptyset$ and, by~\eqref{1.2}, $q\ge 4$. So,
\[
|V(G)|=F_v(2_r;q)\ge F_v(2_r;q-1).
\]
Since $q<r+3$ this contradicts Corollary~\ref{cor:2.1}~(a).
\end{proof}

The following obvious equalities
\begin{align}\label{2.1}
\chi(G_1+G_2)&=\chi(G_1)+\chi(G_2);\\
\cl(G_1+G_2)&=\cl(G_1)+\cl(G_2)
\label{2.2}
\end{align}
are used to prove the following Lemma~\ref{l:2.3}.

Let $f(G)=\chi(G)-\cl(G)$. Then it easily follows from~\eqref{2.1}
and~\eqref{2.2} that
\begin{equation}\label{2.3}
f(G_1+G_2)=f(G_1)+f(G_2).
\end{equation}

\begin{lemma}\label{l:2.3}
Let $m$ and $k$ be positive integers such that $m\ge k+3$ and $2m-1<R(m-k,3)$.
Let
\[
F_v(2_r;r-k+1)\ge r+m
\text{ for any }
r\ge m-1.
\]
Then
\[
F_v(2_r;r-k+1)=r+m
\text{ if }
r\ge m-1.
\]
\end{lemma}

\begin{rem}\label{rem:2.1}
It follows from $r\ge m-1$ and $m\ge k+3$ that $r-k+1\ge 3$.
Thus, by~\eqref{1.2}, the number $F_v(2_r;r-k+1)$ exists.
\end{rem}

\begin{proof}
We need to prove that
\[
F_v(2_r;r-k+1)\le r+m
\text{ if }
r\ge m-1.
\]

It follows from $0<2m-1<R(m-k,3)$ that there exists a graph $P$ such that
$|V(P)|=2m-1$, $\cl(P)\le m-k-1$ and $\alpha(G)<3$. Define
\[
P(r)=K_{r-m+1}+P,
\quad
r\ge m-1.
\]
Since $|V(P)|=2m-1$ and $\alpha(P)<3$ we have $\chi(P)\ge m$. From~\eqref{2.1}
we see that $\chi(P(r))\ge r+1$. The inequality $\cl(G)\le m-k-1$ together
with~\eqref{2.2} implies that $\cl(P(r))\le r-k$. Hence, by~\eqref{1.3},
$P(r)\in H_v(2_r;r-k+1)$ and
\[
F_v(2_r;r-k+1)\le|V(P(r))|=r+m
\text{ if }
r\ge m-1.
\]

Lemma~\ref{2.3} is proved.
\end{proof}

\begin{rem}\label{rem:2.2}
It is clear from the proof of Lemma~\ref{l:2.3} that the following theorem
is true:

\begin{thm}\label{th:2.1}
Let $m$ and $k$ be positive integers such that
\[
2m-1<R(m-k,3)
\text{ and }
m\ge k+3.
\]
Then $F_v(2_r;r-k+1)\le r+m$ if $r\ge m-1$.
\end{thm}
\end{rem}

\section{Some Properties of the Minimal Graphs in $M(x,y)$}\label{s:3}

Let $x$ and $y$ be integers. Define
\[
M(x,y)=\{G:|V(G)|<\chi(G)+2f(G)-x\text{ and }f(G)\le y\}.
\]

In this section we shall prove some properties of the minimal graphs in $M(x,y)$
(see Definition~\ref{d:1.1}). These properties will be need for the proofs of
Theorem~\ref{th:4.1} and Theorem~\ref{th:4.2} in the Section~\ref{s:4}.
If $x<0$ then the empty graph belongs to $M(x,y)$ and hence it is the only
minimal graph in $M(x,y)$. That is why we shall assume $x\ge 0$.

The aim of this section is to prove the following result:

\begin{thm}\label{th:3.1}
Let $M(x,y)\ne\emptyset$, $x\ge 0$ and let $G_0$ be a minimal graph in $M(x,y)$.
If $G_0$ is a nonseparable graph then:
\begin{enumerate}[\indent(a)]
\item
$|V(G_0)|=4f(G_0)-2x-1$;
\item
$4f(G_0)-2x-1<R(f(G_0)-x+1,3)$ where $R(p,3)$ is the Ramsey number.
\end{enumerate}
\end{thm}

An important result of Gallai that we shall need later is:

\begin{thm}[\cite{Gal63} (see also \cite{Gal64})]\label{th:3.2}
Let $G$ be a vertex-critical chromatic graph and $\chi(G)\ge 2$.
Then, if $|V(G)|<2\chi(G)-1$, the graph $G$ is separable in the sense
that $G=G_1+G_2$, where $V(G_i)\ne\emptyset$, $i=1,2$.
\end{thm}

\begin{rem}\label{rem:3.1}
In the original statement of Theorem~\ref{th:3.2} the graph $G$ is
edge-critical (and not vertex-critical) chromatic. Since each vertex-critical
chromatic graph $G$ contains an edge-critical chromatic subgraph $H$ such that
$\chi(H)=\chi(G)$ and $V(H)=V(G)$ the above statement of Theorem~\ref{th:3.2}
is equivalent to the original one.
\end{rem}

In the proof of Theorem~\ref{th:3.1} we shall use the following two Lemmas.

\begin{lemma}\label{l:3.1}
Let $M(x,y)\ne\emptyset$, $x\ge 0$ and $G_0$ be a minimal graph in $M(x,y)$.
Let $A\ne\emptyset$ be an independent vertex set of $G_0$ and $G'_0=G_0-A$.
Then
\begin{enumerate}[\indent\rm(a)]
\item
$\chi(G'_0)=\chi(G_0)-1$;
\item
$\cl(G'_0)=\cl(G_0)$;
\item
$|V(G_0)|=\chi(G_0)+2f(G_0)-x-1$.
\end{enumerate}
\end{lemma}

\begin{proof}
\emph{Proof of} ($a$).
Since $A$ is an independent vertex set we have $\chi(G'_0)=\chi(G_0)-1$ or
$\chi(G'_0)=\chi(G_0)$. Assume that (a) is wrong. Then $\chi(G'_0)=\chi(G_0)$.
Let $\chi(G'_0)=\chi(G_0)=m$ and
\[
V(G'_0)=V_1\cup\dots V_m,
\qquad
V_i\cap V_j=\emptyset,
\quad
i\ne j,
\]
where $V_i$ are independent sets of $G_0$. Note that $\cl(G'_0)\le\cl(G)\le m$.
Thus, after adding new edges $[u,v]$, where $u$ and $v$ belong to different
sets $V_i$ and $V_j$ to $E(G'_0)$ we shall obtain the graph $G''_0$ such
that $\cl(G''_0)=\cl(G_0)$, $\chi(G''_0)=\chi(G_0)$ and $f(G''_0)=f(G_0)$.
Since $A\ne\emptyset$ we have
\[
|V(G''_0)|<|V(G_0)|<\chi(G_0)+2f(G_0)-x=\chi(G''_0)+2f(G''_0)-x.
\]
So, we obtain that $G''_0\in M(x,y)$. This contradicts the minimality of
$G_0$ in $M(x,y)$.

\emph{Proof of} ($b$).
It is clear that $\cl(G'_0)=\cl(G)$ or $\cl(G'_0)=\cl(G_0)-1$. Assume that
(b) is wrong. Then $\cl(G'_0)=\cl(G_0)-1$. By (a) we have $\chi(G'_0)=
\chi(G_0)-1$. Thus, $f(G'_0)=f(G_0)\le y$. Since $|V(G'_0)|<|V(G_0)|$,
from the minimality of $G_0$ it follows that
\[
|V(G'_0)|\ge\chi(G'_0)+2f(G'_0)-x=\chi(G_0)-1+2f(G_0)-x.
\]
From this inequality it follows that $|V(G_0)|\ge\chi(G_0)+2f(G_0)-x$.
This is a contradiction because $G_0\in M(x,y)$.

\emph{Proof of} ($c$).
Assume the opposite, i.e.,
\begin{equation}\label{3.1}
|V(G_0)|\le\chi(G_0)+2f(G_0)-x-2.
\end{equation}
Since $|V(G_0)|\ge\chi(G_0)$ and $x\ge 0$ it follows from~\eqref{3.1} that
$f(G_0)\ne 0$. Thus, there are two non-adjacent vertices $u,v\in V(G_0)$.
Consider the subgraph $G'_0=G_0-\{u,v\}$. By (a) and (b) we have
$\chi(G'_0)=\chi(G_0)-1$ and $f(G'_0)=f(G_0)-1$. Since $|V(G'_0)|=
V(G_0)-2$, it is easy to see from~\eqref{3.1} that
\[
|V(G'_0)|\le\chi(G_0)-1+2f(G_0)-2-x-1<\chi(G'_0)+2f(G'_0)-x.
\]
This is a contradiction since $|V(G'_0)|<|V(G_0)|$.
\end{proof}

\begin{lemma}\label{l:3.2}
Let $M(x,y)\ne\emptyset$, $x\ge 0$ and let $G_0$ be minimal graph in $M(x,y)$.
Then
\begin{enumerate}[\indent\rm(a)]
\item
$G_0$ is a $(\cl(G_0)+1,3)$-graph;
\item
$|V(G_0)|\le 2\chi(G_0)-1$;
\item
$|V(G_0)|\ge 4f(G_0)-2x-1$.
\end{enumerate}
\end{lemma}

\begin{proof}
\emph{Proof of} ($a$).
We need to prove that $\alpha(G_0)<3$. Assume the opposite and let $\{u,v,w\}$
be an independent vertex set of $G_0$. Consider the subgraph $G'_0=G_0-
\{u,v,w\}$. By Lemma~\ref{l:3.1}, we have $\chi(G'_0)=\chi(G_0)-1$ and
$f(G'_0)=f(G_0)-1$. Since $f(G'_0)<y$ and $|V(G'_0)|<|V(G_0)|$, it follows
from the minimality of $G_0$ that
\[
|V(G'_0)|\ge\chi(G'_0)+2f(G'_0)-x.
\]
As $|V(G_0)|=|V(G'_0)|+3$ it follows that $|V(G_0)|\ge\chi(G_0)+2f(G_0)-x$.
This contradicts $G_0\in M(x,y)$.

\emph{Proof of} ($b$).
By (a), $\alpha(G_0)<3$. Thus, we have $|V(G_0)|\le 2\chi(G_0)$ and we need
to prove that $|V(G_0)|\ne 2\chi(G_0)$. Assume the opposite, i.e.,
$|V(G_0)|=2\chi(G_0)$ and let $v\in V(G_0)$. Consider the subgraph $G'_0=
G_0-v$. By Lemma~\ref{l:3.1}~(a), $\chi(G'_0)=\chi(G_0)-1$. Since
$\alpha(G'_0)<3$ it follows that $|V(G'_0)|\le 2\chi(G'_0)-2$ which is
a contradiction.

\emph{Proof of} ($c$).
From (b) and Lemma~\ref{l:3.1}~(c) we obtain
\[
\chi(G_0)\ge 2f(G_0)-x.
\]
By this inequality and Lemma~\ref{l:3.1}~(c) we see that
\[
|V(G_0)|\ge 4f(G_0)-2x-1.\qedhere
\]
\end{proof}

\begin{proof}[\bf Proof of Theorem~\ref{th:3.1}]
\emph{Proof of} ($a$).
According to Lemma~\ref{l:3.1}~(a) $G_0$ is a vertex-critical chromatic graph.
Since $G_0$ is nonseparable, it follows from Lemma~\ref{l:3.2}~(b) and
Theorem~\ref{th:3.2} that
\begin{equation}\label{3.2}
|V(G_0)|=2\chi(G_0)-1.
\end{equation}
By~\eqref{3.2} and Lemma~\ref{l:3.1}~(c) we obtain
\begin{equation}\label{3.3}
\chi(G_0)=2f(G_0)-x,
\quad
\cl(G_0)=f(G_0)-x
\quad\text{and}\quad
|V(G_0)|=4f(G_0)-2x-1.
\end{equation}

\emph{Proof of} ($b$).
According to Lemma~\ref{l:3.2}~(a) we have
\[
|V(G_0)|<R(\cl(G_0)+1,3).
\]
From this inequality and~\eqref{3.3} it follows (b).

Theorem~\ref{th:3.1} is proved.
\end{proof}

\section{A Lower Bound for $|V(G)|$ when $f(G)\le 13$}\label{s:4}

In this section our goal is to prove the following two theorems.

\begin{thm}\label{th:4.1}
Let $G$ be a graph such that $f(G)\le 11$. Then
\begin{enumerate}[\indent\rm(a)]
\item
$|V(G)|\ge\chi(G)+2f(G)$ if $f(G)\le 6$;
\item
$|V(G)|\ge\chi(G)+2f(G)-1$ if $f(G)=7$ or $f(G)=8$;
\item
$|V(G)|\ge\chi(G)+16$ if $f(G)=9$;
\item
$|V(G)|\ge\chi(G)+2f(G)-3$ if $f(G)=10$ or $f(G)=11$.
\end{enumerate}
\end{thm}

\begin{thm}\label{th:4.2}
Let $G$ be a graph such that $f(G)\le 13$. Then
\begin{enumerate}[\indent\rm(a)]
\item
$|V(G)|\ge\chi(G)+2f(G)-4$;
\item
If $f(G)=12$ and $R(10,3)\le 41$ then the inequality \textup{(a)}
is strict.
\end{enumerate}
\end{thm}

\begin{rem}\label{rem:4.1}
If $f(G)\ge 7$ then the inequality (a) of Theorem~\ref{th:4.1} is not true.
For example if $G$ is a minimal graph in $H_v(2_r;r-5)$ we have from
Lemma~\ref{l:2.2} that $\chi(G)=r+1$, $\cl(G)=r-6$ and $f(G)=7$.
By Theorem~\ref{th:1.2} we see that
\[
|V(G)|=r+14<\chi(G)+2f(G)
\text{ if }
r\ge 13.
\]
In the same way we also see that the conditions for $f(G)$ in the statements
(b), (c) and (d) of Theorem~\ref{th:4.1} are necessary.
\end{rem}

\begin{rem}\label{rem:4.2}
If $f(G)\le 6$ the inequality (a) of Theorem~\ref{th:4.1} is exact. Indeed,
if $G$ is a minimal graph in $H_v(2_r;r-k+1)$ where $-1\le k\le 5$, by
Lemma~\ref{l:2.2} we have $\chi(G)=r+1$, $\cl(G)=r-k$ and $f(G)=k+1\le 6$.
When $r$ is large enough we have according to Theorem~\ref{1.1}
\[
|V(G)|=r+2k+3=\chi(G)+2f(G).
\]
In the same way (using Theorem~\ref{th:1.2}) we see that the inequalities
(b), (c) and (d) are exact.
\end{rem}

\begin{rem}\label{rem:4.3}
If $f(G)=13$ the inequality~(a) of Theorem~\ref{th:4.2} is exact by
Theorem~\ref{th:1.4}~(b). If $f(G)=12$ and $R(10,3)\ge 42$ this
inequality is exact according to Theorem~\ref{th:1.3}~(a).
\end{rem}

We shall use the following two lemmas in the proof of Theorem~\ref{th:4.1}
and Theorem~\ref{th:4.2}.

\begin{lemma}\label{l:4.1}
Let $M(0,y)\ne\emptyset$. Then every minimal graph in $M(0,y)$ is nonseparable.
\end{lemma}

\begin{proof}
Assume the opposite and let $G_0$ be a minimal graph in $M(0,y)$ such that
$G_0=G_1+G_2$, where $V(G_i)\ne\emptyset$, $i=1,2$. Since $|V(G_i)|<|V(G_0)|$
we have $G_i\notin M(0,y)$. Since $f(G_i)\le f(G)\le y$ it follows that
\[
|V(G_i)|\ge\chi(G_i)+2f(G_i),
\quad
i=1,2.
\]
Summing these two inequalities we obtain, by~\eqref{2.1} and~\eqref{2.3}, that
\[
|V(G_0)|\ge\chi(G_0)+2f(G_0)
\]
a contradiction.
\end{proof}

\begin{cor}\label{cor:4.1}
$M(0,y)=\emptyset$ if $y\le 6$.
\end{cor}

\begin{proof}
Assume the opposite, i.e., $M(0,y)\ne\emptyset$ for some $y\le 6$. Let $G_0$
be minimal in $M(0,y)$. Then $f(G_0)\le 6$. According to Lemma~\ref{l:4.1}
$G_0$ is nonseparable. Thus, by Theorem~\ref{3.1}~(b) ($x=0$) we have
\[
4f(G_0)-1<R(f(G_0)+1,3)
\]
for
$f(G_0)\le 6$ which is a contradiction (see Table~\ref{t:1.1}).
\end{proof}

\begin{cor}\label{cor:4.2}
Let $G$ be a graph such that
\[
|V(G)|<\chi(G)+2f(G).
\]
Then $|V(G)|\ge 27$.
\end{cor}

\begin{proof}
Since $G\in M(0,f(G))$ we have $M(0,f(G))\ne\emptyset$. Let $G_0$ be a minimal
graph in $M(0,f(G))$. By Corollary~\ref{cor:4.1}, $f(G_0)\ge 7$.
Thus, it follows from Lemma~\ref{l:3.2}~(c) that $V(G)|\ge|V(G_0)|\ge 27$.
\end{proof}

\begin{lemma}\label{l:4.2}
Let $M(x,y)\ne\emptyset$ where $x\ge 0$ and $y\le 13$. Then every minimal graph
in $M(x,y)$ is nonseparable.
\end{lemma}

\begin{proof}
Assume the opposite and let $G_0$ be a minimal graph in $M(x,y)$ such that
$G_0=G_1+G_2$, $V(G_i)\ne\emptyset$, $i=1,2$. Let $f(G_1)\le f(G_2)$. Then
$f(G_1)\le 6$ because $f(G_1)+f(G_2)=f(G_0)\le 13$. By Corollary~\ref{cor:4.1}
we obtain that
\begin{equation}\label{4.1}
|V(G_1)|\ge\chi(G_1)+2f(G_1).
\end{equation}
Since $G_2\notin M(x,y)$ and $f(G_2)\le y$ we have that
\begin{equation}\label{4.2}
|V(G_2)|\ge\chi(G_2)+2f(G_2)-x.
\end{equation}
Summing the inequalities~\eqref{4.1} and~\eqref{4.2} we obtain
by~\eqref{2.1} and~\eqref{2.3} that
\[
|V(G_0)|\ge\chi(G_0)+2f(G_0)-x,
\]
which is a contradiction.
\end{proof}

\begin{proof}[\bf Proof of Theorem~\ref{th:4.1}]
Statement~(a) follows immediately from Corollary~\ref{cor:4.1}.

\emph{Proof of}~(b).
Assume the opposite. Then $M(1,8)\ne\emptyset$. Let $G_0$ be a minimal graph
in $M(1,8)$. It is easy to see that
\[
G_0\in M(1,8)
\Rightarrow
G_0\in M(0,8).
\]
Thus, by Corollary~\ref{cor:4.1}, we have $f(G_0)\ge 7$, i.e.,
$f(G_0)=7$ or $f(G_0)=8$. According to Lemma~\ref{l:4.2} $G_0$
is nonseparable. Thus, from Theorem~\ref{th:3.1} ($x=1$),
it follows that
\[
4f(G_0)-3<R(f(G_0),3),
\]
where $f(G_0)=7$ or $f(G_0)=8$, which is a contradiction.

The proofs of statements (c) and (d) are completely similar to that of
statement~(b).

Theorem~\ref{th:4.1} is proved.
\end{proof}

\begin{proof}[\bf Proof of Theorem~\ref{th:4.2}]
\emph{Proof of}~(a).
Assume the opposite. Then $M(4,13)\ne\emptyset$. Let $G_0$ be a minimal graph
in $M(4,13)$. It is clear that
\[
G_0\in M(4,13)
\Rightarrow
G_0\in M(3,13).
\]
Thus, it follows from Theorem~\ref{4.1} that $f(G_0)\ge 12$. Hence $f(G_0)=12$
or $f(G_0)=13$. By Lemma~\ref{l:4.2}, $G_0$ is nonseparable. Thus,
Theorem~\ref{th:3.1}~(b) ($x=4$) implies
\[
4f(G_0)-9<R(f(G_0)-3,3),
\]
where $f(G_0)=12$ or $f(G_0)=13$ which is a contradiction.

\emph{Proof of}~(b).
Assume the opposite. Then $M(3,12)\ne\emptyset$. Let $G_0$ be a minimal graph
in $M(3,12)$. From Theorem~\ref{th:4.1} it follows that $f(G_0)=12$. Since
$G_0$, by Lemma~\ref{l:4.2}, is nonseparable it follows from
Theorem~\ref{th:3.1}~(b) that
\[
4f(G_0)-7<R(f(G_0)-2,3),
\]
where $f(G_0)=12$ which is a contradiction, by our assumption $R(10,3)\le 41$.
\end{proof}

\section{Proof of Theorem~\ref{th:1.1} and Theorem~\ref{th:1.2}}\label{s:5}

\begin{proof}[\bf Proof of Theorem~\ref{th:1.1}]
\emph{Proof of}~(a).
Let $G$ be a minimal in $H_v(2_r;r-k+1)$. By Lemma~\ref{l:2.2} $\chi(G)=r+1$,
$\cl(G)=r-k$ and $f(G)=k+1$. Since $k\le 5$ we have $f(G)\le 6$. Thus, from
Theorem~\ref{th:4.1}~(a) it follows that
\[
F_v(2_r;r-k+1)=|V(G)|\ge r+2k+3.
\]

\emph{Proof of}~(b).
We shall consider the following three cases.

\textsc{Case~1.}
$k=-1$.
In this case (b) follows from~\eqref{1.4}.

\textsc{Case~2.}
$k\in\{0,2,3,4,5\}$.
By Table~\ref{t:1.1} in this case the following inequality
\[
2(2k+3)-1<R(k+3,3).
\]
holds. Thus, by Lemma~\ref{l:2.3} we obtain $F_v(2_r;r-k+1)=r+2k+3$
if $r\ge 2k+2$.

\textsc{Case~3.}
$k=1$.
We need to prove that $F_v(2_r;r)\le r+5$ if $r\ge 5$. Define
\[
P(r)=K_{r-5}+C_5+C_5,
\quad
r\ge 5.
\]
By~\eqref{2.1} and~\eqref{2.2} we have $\chi(P(r))=r+1$ and $\cl(P(r))=r-1$.
Thus, from~\eqref{1.3} it follows that $P(r)\in H_v(2_r;r)$. Hence
\[
F_v(2_r;r)\le|V(P(r))|=r+5,
\quad
r\ge 5
\]
and Theorem~\ref{th:1.1} is proved.
\end{proof}

\begin{proof}[\bf Proof of Theorem~\ref{th:1.2}]
\emph{Proof of}~(a).
Let $G$ be a minimal graph in $H_v(2_r;r-5)$. Then, by Lemma~\ref{l:2.2},
$\chi(G)=r+1$, $\cl(G)=r-6$ and $f(G)=7$. Thus, from Theorem~\ref{th:4.1}~(b)
it follows
\[
F_v(2_r;r-5)=|V(G)|\ge r+14.
\]
Applying Lemma~\ref{l:2.3} ($k=6$, $m=14$) we obtain
\[
F_v(2_r;r-5)=r+14
\text{ if }
r\ge 13.
\]

Let $8\le r\le 12$. From Table~\ref{t:1.1} we see that $R(r-5,3)\le r+14$.
By Theorem~\ref{th:1.1} ($k=5$) we have $F_v(2_r;r-4)\ge r+13$ and thus
$F_v(2_r;r-4)+1\ge R(r-5,3)$. According to Corollary~\ref{cor:2.1}~(b)
($q=r-4$), $F_v(2_r;r-5)\ge r+15$.

\emph{Proof of}~(b).
Let $G$ be a minimal graph in $H_v(2_r;r-6)$. By Lemma~\ref{l:2.2},
$\chi(G)=r+1$ and $f(G)=8$. From Theorem~\ref{th:4.1}~(b) it follows that
\[
F_v(2_r;r-6)=|V(G)|\ge r+16.
\]
Thus, Lemma~\ref{l:2.3} ($k=7$, $m=16$) implies $F_v(2_r;r-6)=r+16$ if
$r\ge 15$.

\emph{Proof of}~(c).
Let $G$ be a minimal graph in $H_v(2_r;r-7)$. By Lemma~\ref{l:2.2},
$\chi(G)=r+1$ and $f(G)=9$. Thus, from Theorem~\ref{th:4.1}~(c)
it follows that $F_v(2_r;r-7)\ge r+17$, $r\ge 10$. From this inequality
and Lemma~\ref{l:2.3} ($k=8$, $m=17$) we see that $F_v(2_r;r-7)=r+17$
if $r\ge 16$.

Let $10\le r\le 15$. By Table~\ref{t:1.1} we have that $R(r-7,3)<r+17$.
Since, by~(b), $F_v(2_r;r-6)+1\ge r+17$ we have $F_v(2_r;r-6)+1>R(r-7,3)$.
From Corollary~\ref{cor:2.1}~(b), the inequality $F_v(2_r;r-7)\ge r+18$
holds.

\emph{Proof of}~(d).
If $G$ be a minimal graph in $H_v(2_r;r-8)$ then, by Lemma~\ref{l:2.2},
$\chi(G)=r+1$ and $f(G)=10$. From Theorem~\ref{th:4.1}~(d) it follows
that
\[
|V(G)|=F_v(2_r;r-8)\ge r+18,
\quad
r\ge 11.
\]
Applying Lemma~\ref{l:2.3} ($k=9$, $m=18$) we obtain $F_v(2_r;r-8)=r+18$
if $r\ge 17$.

Let $11\le r\le16$. In this case we have $R(r-8,3)\le r+18$. By~(c),
$F_v(2_r;r-7)\ge r+17$. Thus, $F_v(2_r;r-7)+1\ge R(r-8,3)$ and,
by Corollary~\ref{cor:2.1}~(b), $F_v(2_r;r-8)\ge r+19$.

\emph{Proof of}~(e).
Let $G$ be a minimal graph in $H_v(2_r;r-9)$. According to Lemma~\ref{l:2.2}
we have $\chi(G)=r+1$ and $f(G)=11$. By Theorem~\ref{th:4.1}~(d) we obtain
\[
|V(G)|=F_v(2_r;r-9)\ge r+20.
\]
This inequality and Lemma~\ref{l:2.3} ($k=10$, $m=20$) imply that
$F_v(2_r;r-9)=r+20$ if $r\ge 19$.

Theorem~\ref{th:1.2} is proved.
\end{proof}

\section{Proof of Theorem~\ref{th:1.3} and Theorem~\ref{th:1.4}}\label{s:6}

\begin{proof}[\bf Proof of Theorem~\ref{th:1.3}]
Let $G$ be a minimal graph in $H_v(2_r;r-10)$. According to Lemma~\ref{l:2.2}
we have $\chi(G)=r+1$ and $f(G)=12$. Thus, by Theorem~\ref{th:4.2}~(a) it
follows that
\[
|V(G)|=F_v(2_r;r-10)\ge r+21,
\quad
r\ge 13.
\]

Let $R(10,3)>41$. Then, by Lemma~\ref{l:2.3} ($k=11$, $m=21$) it follows that
\[
F_v(2_r;r-10)=r+21
\text{ if }
r\ge 20.
\]

Let $R(10,3)\le 41$. From Theorem~\ref{th:4.2}~(b) we obtain
$|V(G)|=F_v(2_r;r-10)\ge r+22$. Applying Lemma~\ref{l:2.3}
($k=11$, $m=22$) we deduce that $F_v(2_r;r-10)=r+22$ if $r\ge 21$
because $43<R(11,3)$ (see~\cite{Rad06}).
\end{proof}

\begin{proof}[\bf Proof of Theorem~\ref{th:1.4}]
\emph{Proof of}~(a).
The proof is by induction on $k$ with induction base $k=12$. Let $G$ be
a minimal graph in $H_v(2_r;r-11)$. Then, by Theorem~\ref{l:4.2}~(a) we obtain
\begin{equation}\label{6.1}
|V(G)|=F_v(2_r;r-11)\ge r+23.
\end{equation}
We are done with the base $k=12$. Let $k\ge 13$ and
\[
F_v(2_r;r-k+2)\ge r+k+10.
\]
Then, by Corollary~\ref{cor:2.1}~(a) it follows that
\[
F_v(2_r;r-k+1)\ge r+k+11.
\]

\emph{Proof of}~(b).
From~\eqref{6.1} and Lemma~\ref{l:2.3} ($k=12$, $m=23$) we deduce that
$F_v(2_r;r-11)=r+23$ if $r\ge 22$ because $R(11,3)>45$ (see~\cite{Rad06}).

Theorem~\ref{th:1.4} is proved.
\end{proof}

\section{Lower Bounds for Arbitrary Vertex Folkman numbers}\label{s:7}

Let $a_1,\dots,a_r$ be positive integers. Define
\begin{equation}\label{7.1}
m(a_1,\dots,a_r)=m=\sum_{i=1}^r(a_i-1)+1.
\end{equation}
It is easy to see that $K_m\tov(a_1,\dots,a_r)$ and
$K_{m-1}\notov(a_1,\dots,a_r)$. Therefore
\[
\FN=m
\text{ if }
q>m.
\]
By~\eqref{1.1}. the Folkman number $\FN[m]$ exists only when
$m\ge\max\{a_1,\dots,a_r\}+1$. It was proved in~\cite{LucRucUrb}
that
\[
\FN[m]=m+\max\{a_1,\dots,a_r\}.
\]
The exact values of all numbers $\FN[m-1]$ for which
$\max\{a_1,\dots,a_r\}\le 4$ are known. A detailed
exposition of these results was given in~\cite{Nen02}.
We must add the equality $F_v(2,2,3;4)=14$ obtained
in~\cite{ColRad06}. We do not know any exact values of
$\FN[m-1]$ in the case when $\max\{a_1,\dots,a_r\}\ge 5$.

In this section we shall use the following result \cite{Nen01}
\begin{equation}\label{7.2}
G\tov(a_1,\dots,a_r)
\Rightarrow
\chi(G)\ge m.
\end{equation}

Let $G$ be a minimal graph in $H_v(a_1,\dots,a_r;q)$. Then, be~\eqref{7.2}
and~\eqref{1.3} it follows that $G\in H_v(2_{m-1};q)$. Thus we have
$|V(G)|\ge F_v(2_{m-1};q)$. So, we obtain
\begin{equation}\label{7.3}
\FN\ge F_v(2_{m-1};q),
\end{equation}
where $m$ is defined by the equality~\eqref{7.1}. From~\eqref{7.3},
Theorem~\ref{th:1.1}, Theorem~\ref{th:1.2}, Theorem~\ref{th:1.3} and
Theorem~\ref{th:1.4} we easily get the following theorem:

\begin{thm}\label{th:7.1}
Let $a_1,\dots,a_r$ be integers, $a_i\ge 2$, $i=1,\dots,r$
and $m=\sum_{i=1}^r(a_i-1)+1$. Let $k$ be an integer such that
\begin{equation}\label{7.4}
m-k>\max\{a_1,\dots,a_r\}.
\end{equation}
Then the following inequalities hold:
\begin{align*}
&\FN[m-k]\ge m+2k+2\text{ if }-1\le k\le 5;\\
&\FN[m-6]\ge m+13;\\
&\FN[m-7]\ge m+15;\\
&\FN[m-8]\ge m+16;\\
&\FN[m-9]\ge m+17;\\
&\FN[m-10]\ge m+19;\\
&\FN[m-11]\ge m+20;\\
&\FN[m-11]\ge m+21\text{ if }R(10,3)\le 41;\\
&\FN[m-k]\ge m+k+10\text{ if }k\ge 12.
\end{align*}
\end{thm}

\begin{rem}\label{rem:7.1}
According to~\eqref{1.1} the inequality~\eqref{7.4} in the statement of
Theorem~\ref{th:7.1} is necessary.
\end{rem}

\begin{proof}
Since all inequalities are proved in the same way, we shall prove the last one
only. By Theorem~\ref{th:1.4} we have
\begin{equation}\label{7.5}
F_v(2_r;r-k+1)\ge r+k+11,
\quad
r\ge k+2.
\end{equation}
As $\max\{a_1,\dots,a_r\}\ge 2$, it follows from~\eqref{7.4} that $m-1\ge k+2$.
Thus, the inequality~\eqref{7.5} is true for $r=m-1$, i.e.,
\begin{equation}\label{7.6}
F_v(2_{m-1};m-k)\ge m+k+10.
\end{equation}
We obtain from~\eqref{7.6} and~\eqref{7.3} that
\[
F_v(a_1,\dots,a_r;m-k)\ge m+k+10.\qedhere
\]
\end{proof}

\begin{rem}\label{rem:7.2}
Dudek and R\"odl \cite{Dudek08-4} proved that
\[
\FN\le cp^3\log^3p,
\]
where $p=\max\{a_1,\dots,a_r\}$ and $c$ is a constant depending only on $r$.
\end{rem}

\section{Lower Bounds for Edge Folkman Numbers}\label{s:8}

Let $a_1\dots,a_r$ be integers, $a_i\ge 2$. The symbol $G\tov[e](a_1,\dots,a_r)$
denotes that in every $r$-coloring of the edge set $E(G)$ there exists
a monochromatic $a_i$-clique of color $i$ for some $i\in\{1,\dots,r\}$. Define
\begin{align*}
H_e(a_1,\dots,a_r;q)&=\{G:G\tov[e](a_1,\dots,a_r)\text{ and }\cl(G)<q\},\\
F_e(a_1,\dots,a_r;q)&=\min\{|V(G)|:G\in H_e(a_1,\dots,a_r;q)\}.
\end{align*}

It is clear that from $G\tov[e](a_1,\dots,a_r)$ it follows $\cl(G)\ge
\max\{a_1,\dots,a_r\}$. There exists a graph $G\tov[e](a_1,\dots,a_r)$
and $\cl(G)=\max\{a_1,\dots,a_r\}$. In the case $r=2$ this was proved
in~\cite{Folkman} and the general case in~\cite{Nes76}. Thus, we have
\begin{equation}\label{8.1}
F_e(a_1,\dots,a_r;q)\text{ exists}
\iff
q>\max\{a_1,\dots,a_r\}.
\end{equation}
The numbers $F_e(a_1,\dots,a_r;q)$ are called \emph{edge Folkman numbers}.

From definition of Ramsey number $R(a_1,\dots,a_r)$ it follows that
\[
F_e(a_1,\dots,a_r;q)=R(a_1,\dots,a_r)
\text{ if }
q>R(a_1,\dots,a_r).
\]
Thus, we consider only numbers $F_e(a_1,\dots,a_r;R(a_1,\dots,a_r)-k)$,
where $k\ge -1$. An exposition of the known edge Folkman numbers is given
in~\cite{KolNen}. We must add the new upper bounds for the number
$F_e(3,3;4)$ obtained in~\cite{Dudek08-5} and~\cite{Lu}.

In this section we shall use the following result obtained
by S.~Lin~\cite{Lin}
\begin{equation}\label{8.2}
G\tov[e](a_1,\dots,a_r)
\Rightarrow
\chi(G)\ge R(a_1,\dots,a_r).
\end{equation}

From~\eqref{8.2} and~\eqref{1.3} we see that
\[
G\in H_e(a_1,\dots,a_r;q)
\Rightarrow
G\in H_v(2_{R-1};q),
\]
where $R=R(a_1,\dots,a_r)$. Thus, we have
\begin{equation}\label{8.3}
F_e(a_1,\dots,a_r;q)\ge F_v(2_{R-1};q).
\end{equation}

From~\eqref{8.3}, Theorem~\ref{th:1.1}, Theorem~\ref{th:1.2},
Theorem~\ref{th:1.3} and Theorem~\ref{th:1.4} it easily follows
the following statement.

\begin{thm}\label{th:8.1}
Let $a_1,\dots,a_r$ be integers, $a_i\ge 2$, $i=1,\dots,r$. Let
\[
R-k>\max\{a_1,\dots,a_r\},
\]
where $k\ge -1$ is integer and $R=R(a_1,\dots,a_r)$. Then
\begin{align*}
&F_e(a_1,\dots,a_r;R-k)\ge R+2k+2\text{ if }-1\le k\le 5;\\
&F_e(a_1,\dots,a_r;R-6)\ge R+13;\\
&F_e(a_1,\dots,a_r;R-7)\ge R+15;\\
&F_e(a_1,\dots,a_r;R-8)\ge R+16;\\
&F_e(a_1,\dots,a_r;R-9)\ge R+17;\\
&F_e(a_1,\dots,a_r;R-10)\ge R+19;\\
&F_e(a_1,\dots,a_r;R-11)\ge R+20;\\
&F_e(a_1,\dots,a_r;R-11)\ge R+21\text{ if }R(10,3)\le 41;\\
&F_e(a_1,\dots,a_r;R-k)\ge R+k+10\text{ if }k\ge 12.
\end{align*}
\end{thm}

\begin{rem}\label{rem:8.1}
According to~\eqref{8.1} the inequality
\[
R-k>\max\{a_1,\dots,a_r\}
\]
in the statement of Theorem~\ref{th:8.1} is necessary.
\end{rem}

\begin{rem}\label{rem:8.2}
In the particular cases $k=0$ and $k=1$ Theorem~\ref{th:8.1} was proved
by S.~Lin \cite{Lin}. Lin~\cite{Lin} also proved that when $k=0$
the respective inequality in Theorem~\ref{th:8.1} is exact and
the conjecture was raised that if $k=1$ the first inequality in
Theorem~\ref{th:8.1} is strict. This Lin's hypothesis was disproved
in~\cite{Nen80}, where the equality $F_e(3,3,3;16)=21$ was established.
The particular cases $k=2$ and $k=3$ of Theorem~\ref{th:8.1} were proved
in~\cite{Nen81b} and~\cite{Nen81a}, respectively. In~\cite{Nen81b}
and~\cite{Nen81a} it was also proved that if $k=2$ and $k=3$ then respective
inequalities of Theorem~\ref{th:8.1} are exact. The other inequalities are
new. We do not know whether these inequalities are exact.
\end{rem}

\bibliographystyle{plain}
\bibliography{NNSerdica}

\begin{flushleft}
Faculty of Mathematics and Informatics\\
St. Kl. Ohridski University of Sofia\\
5, J. Bourchier Blvd.\\
BG-1164 Sofia, Bulgaria\\
e-mail: \texttt{nenov@fmi.uni-sofia.bg}
\end{flushleft}

\end{document}